\documentclass[11pt]{article}
\usepackage[utf8]{inputenc}
\usepackage{subfigure}
\usepackage{tikz}
\usepackage[utf8]{inputenc}
\usepackage[english]{babel}
\usepackage{amsmath}
\usepackage{times}
\usepackage{graphicx}
\usepackage{float}
\usepackage[margin=1in]{geometry}
\usepackage{blindtext}
\usepackage{amssymb}
\usepackage{tikz}
\usepackage{amsthm}
\usepackage{mathtools}
\usepackage{textcmds}

\newtheorem{thm}{Theorem}[section]
\newtheorem{cor}[thm]{Corollary}
\newtheorem{lem}[thm]{Lemma}
\newtheorem{pro}[thm]{Problem}

\newtheorem{claim}{Claim}

\title{\bf The $\alpha$-spectral Tur\'an type problems for graphs\thanks{Research was partially supported by the National
		Nature Science Foundation of China (grant numbers 12331012, 12571375)}}
\date{}
\date{}
\author {Jiadong Wu$^{1}$, \, Yongchun Lu$^{1}$,  Liying Kang$^{1,2}$\thanks{\em Corresponding author. Email address: lykang@shu.edu.cn (L. Kang), 1753381890@qq.com (J. Wu), luyongchun@shu.edu.cn (Y. Lu)} \\
	{\small $^{1}$ Department of Mathematics, Shanghai University,
		Shanghai 200444, P.R. China}\\
	{\small$^{2}$Newtouch Center for Mathematics of Shanghai University,
		Shanghai,  China, 200444}}

\begin{document}
	\maketitle

\begin{abstract}
       For $0 \leq \alpha < 1$, the $\alpha$-spectral radius of a graph $G$ is defined as
       the  largest eigenvalue of $A_{\alpha}(G)=\alpha D(G)+(1-\alpha)A(G)$, where $D(G)$ and $A(G)$ are the diagonal matrix of degrees and adjacency matrix of $G$, respectively.
       A graph is called color-critical if it contains an edge whose deletion reduces its chromatic number.
       The celebrated  Erd\H{o}s-Stone-Simonovits theorem asserts that $ \mathrm{ex}(n,\mathcal{F})=\left(1-\frac{1}{\chi(\mathcal{F})-1}+o(1)\right)\frac{n^2}{2},$
   where $\chi(\mathcal{F})$ is the chromatic number of $\mathcal{F}$.
    Nikiforov and Zheng et al.   established the adjacency spectral and signless Laplacian spectral versions of this theorem, respectively.
    In this paper, we present the   $\alpha$-spectral version of this theorem, which unifies the aforementioned  results. Furthermore, we characterize the   $\alpha$-spectral extremal graphs for color-critical graphs, thereby extending the existing results on adjacency spectral and signless Laplacian spectral extremal graphs for such graphs.

		\bigskip
		
		
	\end{abstract}

\section{Introduction}

Let $G=(V(G),E(G))$ be a simple graph with vertex set $V(G)$ and edge set $E(G)$. As usual, we use $|G|$ and $e(G)$ to represent the number of vertices and edges of a graph $G$, respectively. For $v\in V(G)$, the neighbors of $v$ in $G$ are denoted by $N_{G}(v)$, and $d_{G}(v)=|N_{G}(v)|$ denotes its degree in $G$. If there is no ambiguity, we simplify $N_{G}(v)$ and $d_{G}(v)$ as $N(v)$ and $d(v)$, respectively. We write $\Delta(G)$ for the maximum degree of $G$ and $\delta(G)$ for the minimum degree of $G$. Let $T_{n,r}$ denote the $r$-partite Tur\'an graph of order $n$. The  $k$-vertex independent set is denoted by $I_k$.
For two graphs $G$ and $H$, the disjoint union of $G$ and $H$ is denoted by $G\cup H$. The \textit{join} of $G$ and $H$, denoted by $G\vee H$, is the graph obtained from $G\cup H$ by adding all possible edges between $G$ and $H$.
Denote $S_{n,k}=K_{k}\vee I_{n-k}$ and $S^{+}_{n,k}=K_{k}\vee (K_2 \cup I_{n-k-2})$. Let $M_{k+1}$ be a matching of size $k+1$ and $P_k$ be a path of order $k$. A cycle on $k$ vertex is denoted by $C_k$.
A graph is called \textit{color-critical} if it contains an edge whose deletion reduces its chromatic number.

Let $\mathcal{F}$ be a family of graphs. A graph is called $\mathcal{F}$-free if it does not contain any member of $\mathcal{F}$ as a subgraph. The classic \textit{Tur\'an type problem} is to determine the maximum number of edges over all $\mathcal{F}$-free graphs of order $n$. The \textit{Tur\'an number} of $\mathcal{F}$ is defined as follows:
	\begin{align}
		\operatorname{ex}(n,\mathcal{F})=\max\{e(G)|\ G \text{ is an $n$-vertex }\mathcal{F}\text{-free graph}\}.\nonumber
	\end{align}
Let $\chi(\mathcal{F})=\min_{F\in \mathcal{F}} \chi(F)$.
When  $\mathcal{F}=\{F\}$, we simplify the notation to  $\mathrm{ex}(n,F)$ instead of $\mathrm{ex}(n,\{F\})$. The same simplification applies to all other relevant functions. The \textit{Turán density} of $ \mathcal{F} $ is defined as
\[
\pi(\mathcal{F}) := \lim_{n \to \infty} \frac{\mathrm{ex}(n, \mathcal{F})}{\binom{n}{2}},
\]
for which the existence is guaranteed by an averaging argument due to Katona, Nemetz and Simonovits in \cite{Katona}.

 Let $G$ be a graph with adjacency matrix $A(G)$, and let $D(G)$ be the diagonal matrix of the degrees of $G$.  The \textit{signless Laplacian matrix} of a graph $G$ is defined as $Q(G):=D(G)+A(G)$. In 2017, Nikiforov \cite{N1} proposed to study the convex combinations $A_{\alpha}(G)$ of $A(G)$ and $D(G)$ defined by
 $$ A_{\alpha}(G):=\alpha D(G)+(1-\alpha)A(G), \ \ 0\leq \alpha<1.$$
 Obviously, $A(G)=A_0(G)$ and  $Q(G)=2A_{1/2}(G)$. The largest eigenvalue of $A_{\alpha}(G)$ is called the \textit{$\alpha$-spectral radius} of $G$, denoted by $\lambda_{\alpha}(G)$.
 The $\lambda_{0}(G)$ (resp. $2\lambda_{1/2}(G)$) of $G$ is usually referred to as the spectral radius (resp. signless Laplacian spectral radius) of $G$, denoted by $\lambda(G)$ (resp. $q(G)$).
 In 2017, Nikiforov \cite{N1} proposed the following  problem.

\begin{pro}
    Given a family of graphs $\mathcal{F}$, what is the maximum $\lambda_{\alpha}(G)$ of a graph $G$ of order $n$, with no subgraph isomorphic to any $F\in \mathcal{F}$?
\end{pro}

Let  $ \mathrm{ex}_{\alpha}(n, \mathcal{F})$  denote the maximum   $ \alpha $-spectral radius among all   $ \mathcal{F} $-free graphs of order   $ n $. An   $\mathcal{F}$-free graph of order  $n$  that achieves $ \mathrm{ex}_{\alpha}(n, \mathcal{F}) $   is called an \textit{  $\alpha$-spectral extremal graph} for   $\mathcal{F}$, and we let $ \mathrm{Ex}_{\alpha}(n, \mathcal{F}) $   denote the family of all such graphs. For specific spectral parameters, $\mathrm{ex}_{\lambda}(n,\mathcal{F})$   stands for the maximum spectral radius of   $\mathcal{F}$-free graphs of order    $n$, while $\mathrm{ex}_{q}(n,\mathcal{F})$   denotes the maximum signless Laplacian spectral radius of   $\mathcal{F}$-free graphs on  $n$  vertices.

The study of $\alpha$-spectral Tur\'an problem was initiated by Nikiforov \cite{N1}, who determined the $\alpha$-spectral extremal graphs corresponding to  complete graphs. For matching $M_{k+1}$, Yuan and Shao \cite{Y1} showed that $\mathrm{Ex}_{\alpha}(n,M_{k+1})=S_{n,k}$. Further advancing the field,
 Chen, Liu, and Zhang \cite{C1} proved that $\mathrm{Ex}_{\alpha}(n,P_{2k+2})=S_{n,k}$ and $\mathrm{Ex}_{\alpha}(n,P_{2k+3})=S^{+}_{n,k}$,
shedding light on the   $\alpha$-spectral extremal structures for paths of specific lengths.
 For cycles, Li and Yu \cite{L1} showed that $\mathrm{Ex}_{\alpha}(n,C_{2k+2})=S^{+}_{n,k}$ and $\mathrm{Ex}_{\alpha}(n,\{C_{2k+1}, C_{2k+2}\})=S_{n,k}$.
 Expanding to disjoint cycles,
  Li, Yu, and Zhang \cite{L2} determined the $\alpha$-spectral extremal graphs, they obtained $\mathrm{Ex}_{\alpha}(n,\mathcal{F})=S_{n,2k-1}$, where  $\mathcal{F}$ is the family of all disjoint unions of $k$ cycles.
  A nice result was achieved by Chen et al. \cite{C2}, who proved the   $\alpha$-spectral Erd\H{o}s-S\'os theorem. This theorem asserts that
  every $n$-vertex graph $G$ with $\lambda_{\alpha}(G)\geq \lambda_{\alpha}(S_{n,k})$ contains all trees on $2k+2$ vertices unless $G=S_{n,k}$. Moreover, if $\lambda_{\alpha}(G)\geq \lambda_{\alpha}(S^{+}_{n,k})$, then $G$ contains all trees on $2k+3$ vertices unless $G=S^{+}_{n,k}$.
 Most recently, Byrne, Desai, and Tait \cite{B}  established a unified framework that characterizes the   $\alpha$-spectral extremal graphs for a broad class of graph families.


\subsection{The  spectral  version Erd\H{o}s-Stone-Simonovits theorem }

The famous Erd\H{o}s-Stone-Simonovits theorem gives an asymptotic value of the Tur\'an number for a family of graphs.

\begin{thm} [Erd\H{o}s-Stone-Simonovits \cite{E1, E2}] \label{T6}
   If $\mathcal{F}$ is a family of graphs with  $\chi(\mathcal{F})\geq 2$, then
   $$ \mathrm{ex}(n,\mathcal{F})=\left(1-\frac{1}{\chi(\mathcal{F})-1}+o(1)\right)\frac{n^2}{2}.$$
\end{thm}
\noindent
As a direct consequence of Theorem \ref{T6}, we have   $\pi(\mathcal{F})=1-\frac{1}{\chi(\mathcal{F})-1}$. Nikiforov \cite{N2} later proposed a spectral analog of the Erd\H{o}s-Stone-Simonovits theorem, stated as follows.
\begin{thm} [Nikiforov \cite{N2}] \label{T7}
    If $F$ is a graph with $\chi(F)\geq 2$, then
     $$ \mathrm{ex}_{\lambda}(n,F)=\left(1-\frac{1}{\chi(F)-1}+o(1)\right)n.$$
\end{thm}

More recently, Zheng, Li, and Su \cite{Z1} established the signless Laplacian spectral counterpart of this theorem.
\begin{thm} [Zheng, Li, and Su \cite{Z1}] \label{T8}
    If $F$ is a graph with $\chi(F)\geq 3$, then
     $$ \mathrm{ex}_{q}(n,F)=\left(1-\frac{1}{\chi(F)-1}+o(1)\right)2n.$$
\end{thm}

 For a bipartite graph $F$, Zheng, Li,  and Su \cite{Z1} also derived the following result.
 \begin{thm} [Zheng, Li, and Su \cite{Z1}] \label{T9}
If $F$ is a bipartite graph, i.e., $\chi(F)=2$, then
$$\mathrm{ex}_q(n,F)=
\begin{cases}
\bigl(1+o(1)\bigr)\,n & \text{if }F \text{ is not a star};\\[2mm]
2(k-1) & \text{if }F  \text{ is a star } K_{1, k}.
\end{cases}$$
\end{thm}

\subsection{The spectral extremal results for color-critical graphs }


In 2007, Nikiforov \cite{N4} showed that if $G$ is a $K_{r+1}$-free graph on $n$ vertices,
then $\lambda (G)\le \lambda (T_{n,r})$, with equality if and only if
$G=T_{n,r}$.
 Subsequently, He, Jin, and Zhang \cite{H1} extended this line of research to the signless Laplacian spectral radius, establishing the following theorem:
\begin{thm} [He, Jin, and Zhang \cite{H1}]
 If $G$ is  an $n$-vertex $K_{r+1}$-free graph, then
 $$q(G)\leq q(T_{n,r}).$$
Moreover, the equality holds if and only if $G$ is a complete bipartite graph for $r=2$ and the $r$-partite Tur\'an graph $T_{n,r}$ for every $r\geq 3$.
\end{thm}

Later, Nikiforov \cite{N1} generalized the aforementioned results and obtained the following theorem.

\begin{thm} [Nikiforov \cite{N1}] \label{TT0}
 Let $r \geq 2$ and $G$ be a $K_{r+1}$-free graph of order $n$.
\begin{itemize}
    \item[(1)]  If $0 \leq \alpha < 1 - \frac{1}{r}$, then $\lambda_{\alpha}(G) < \lambda_{\alpha}(T_{n,r})$, unless $G = T_{n,r}$.

    \item[(2)]  If $1 > \alpha > 1 - \frac{1}{r}$, then $\lambda_{\alpha}(G)< \lambda_{\alpha}(S_{n,r-1})$, unless $G = S_{n,r-1}$.

    \item[(3)] If $\alpha = 1 - \frac{1}{r}$, then $\lambda_{\alpha}(G) \leq \left(1 - \frac{1}{r}\right)n$, with equality if and only if $G$ is a complete $r$-partite graph.
\end{itemize}
\end{thm}

Color-critical graphs form a broad and fundamental class of graphs. Examples include cliques, odd cycles, book graphs, and even wheels-all of which are color-critical. Let  $F$  be a color-critical graph with   $\chi(F)=r+1$. Simonovits \cite{S1} determined the Tur\'an number of   $F$, proving that $\mathrm{Ex}(n,F)=T_{n,r}$   for sufficiently large   $n$. Nikiforov \cite{N3} established a spectral counterpart to this result, showing that   $\mathrm{ex}_{\lambda}(n,F)=\lambda(T_{n,r})$. Subsequently, Zheng, Li, and Li \cite{Z2} settled the signless Laplacian spectral extremal problem for   $F$, determining  $\mathrm{ex}_{q}(n,F)$  explicitly.

\begin{thm} [Zheng, Li, and Li \cite{Z2}] \label{TZ}
   Let $F$ be a color-critical graph with $\chi(F)=r+1\geq 4$. Then there exists $n_0$ such that for every $F$-free graph $G$ on $n\geq n_0$ vertices, we have
        $$ q(G)\leq q(T_{n,r}), $$
    where the equality holds if and only if $G$ is the $r$-partite Tur\'an graph $T_{n,r}$.
\end{thm}

\section{Main results}

 Given a family of graphs with $\chi(\mathcal{F})=r+1\geq 3$, we first establish the asymptotic value of $\mathrm{ex}_{\alpha}(n, \mathcal{F})$ as follows.
\begin{thm}\label{TT1}
    Let $\mathcal{F}$ be a family of graphs with $\chi(\mathcal{F})=r+1\geq 3$ and $0\leq  \alpha\leq 1- \frac{1}{r}$. Then
    $$\mathrm{ex_{\alpha}}(n, \mathcal{F})=\left(1-\frac{1}{r}+o(1)\right)n.$$
\end{thm}
\noindent\textbf{Remark.} Theorem \ref{TT1} does not hold when $\alpha>1-\frac{1}{r}$. Nikiforov \cite{N1} proved that $\mathrm{Ex}_{\alpha}(n,K_{r+1})=S_{n,r-1}$ for $r\geq 2$ and $1>\alpha>1-\frac{1}{r}$. Clearly, $\lambda_{\alpha}(S_{n,r-1})\geq \alpha(n-1)>\left(1-\frac{1}{r}+o(1)\right)n$ for sufficiently large  $n$.

Suppose $G$ is an $\mathcal{F}$-free graph of order $n$. By Lemma \ref{LL5} (stated in Section 3, which asserts $e(G)\leq  \lambda_{\alpha}(G) \frac{n}{2}$)  and Theorem \ref{TT1}, we have  $e(G)\leq  \lambda_{\alpha}(G) \frac{n}{2} \leq  \left(1-\frac{1}{r}+o(1)\right)\frac{n^2}{2}$. Thus, Theorem \ref{TT1} implies Theorem \ref{T6} in the case where $\chi(\mathcal{F})\geq 3$. For families $\mathcal{F}$ containing a bipartite graph, we establish the following result.

\begin{thm}\label{TT11}
    Let $\mathcal{F}$ be a family of graphs with $\chi(\mathcal{F})=2$ and $0\leq \alpha< 1$. Then
    \begin{itemize}
        \item [(1)] $\mathrm{ex_{\alpha}}(n, \mathcal{F})=\left(\alpha+o(1)\right)n $ if $\mathcal{F}$ does not contain $K_{1,k}\cup s K_1$ for any $k\geq 1, s\geq 0$.
        \item [(2)] $\mathrm{ex_{\alpha}}(n, \mathcal{F})=O(1) $ if $\mathcal{F}$ contains $K_{1,k}\cup s K_1$ for some $k\geq 1, s\geq 0$. Particularly, if $\mathcal{F}=\{K_{1,k}\}$, then $\mathrm{ex_{\alpha}}(n, \mathcal{F})=k-1$.
    \end{itemize}
\end{thm}

  When $\alpha=0$, Theorem \ref{TT1} combined with Theorem \ref{TT11} implies Theorem \ref{T7}. Furthermore, setting   $\alpha=\frac{1}{2}$,  Theorems \ref{TT1} and \ref{TT11} readily yield Theorems \ref{T8} and \ref{T9}.

 For color-critical graphs   $F$, we further determine the maximum   $\alpha$-spectral radius among all   $F$-free graphs of order   $n$.

\begin{thm}\label{TT2}
    Let $F$ be a color-critical graph with $\chi(F)=r+1\geq 3$, and let $G$ be an $F$-free graph of order $n$. For sufficiently large $n$ and $0\leq \alpha< 1-\frac{1}{r}$, we have
    $$\lambda_{\alpha}(G)\leq \lambda_{\alpha}(T_{n,r}),$$
    where   equality holds if and only if $G=T_{n,r}$.
\end{thm}
\noindent\textbf{Remark.} Theorem \ref{TT2} does not  hold when $\alpha=1-\frac{1}{r}$.  Indeed, for this specific value of $\alpha$, Theorem \ref{TT0} asserts that $\mathrm{Ex}_{\alpha}(n,K_{r+1})$ consists of  all complete $r$-partite graphs.

Clearly, Theorem \ref{TT2} extends Theorem \ref{TT0} for $0\leq \alpha< 1-\frac{1}{r}$. When $\alpha=0$, it  recovers Nikiforov's spectral result for color-critical graphs \cite{N3}; when   $\alpha=\frac{1}{2}$, it implies Theorem \ref{TZ}.

The \textit{even wheel} of order $2k+2$ is defined as $W_{2k+2}:=K_{1}\vee C_{2k+1}$. Notably, $W_{2k+2}$ is color-critical with $\chi(W_{2k+2})=4$. By Theorem \ref{TT2},
we immediately obtain the following result.

\begin{cor}
    For $k\geq 1$, $0\leq \alpha<\frac{2}{3}$ and sufficiently large $n$, if $G$ is a $W_{2k+2}$-free graph of order $n$, then
    $$\lambda_{\alpha}(G)\leq \lambda_{\alpha}(T_{n,3}),$$
    where  equality holds  if and only if $G=T_{n,3}$.
\end{cor}

The \textit{generalized book graph} $B_{r,k}$ is defined as $B_{r,k}:=K_{r}\vee I_k$. One can see that $B_{r,k}$ is a color-critical graph and $\chi(B_{r,k})=r+1$. Thus, Theorem \ref{TT2} implies the following result.

\begin{cor}
    Suppose $k\geq 1$ and $r\geq 2$. For sufficiently large $n$ and $0\leq \alpha<1-\frac{1}{r}$, if $G$ is a $B_{r,k}$-free graph of order $n$, then
    $$\lambda_{\alpha}(G)\leq \lambda_{\alpha}(T_{n,r}),$$
    with  equality  if and only if $G=T_{n,r}$.
\end{cor}

\section{Preliminaries}

For a graph $G$ and integer $p\ge 1$, the \textit{blow-up} of $G$, denoted by $G^{p}$, is the graph obtained from $G$ by replacing each vertex of $G$ by a set of $p$ independent vertices and each edge of $G$ by a complete bipartite graph $K_{p,p}$.
A family of graphs $\mathcal{G}$ is  called  \textit{multiplicative} if $G\in \mathcal{G}$ implies $G^{p}\in \mathcal{G}$ for all integers $p\geq 1$.
A family $\mathcal{G}$ is said to be \textit{hereditary} if it is closed under taking induced subgraphs.

Let $\mathcal{F}$ be a family of graphs and $\mathcal{G}$ be a class of $\mathcal{F}$-free graphs. If there exist $\varepsilon>0$ and $N>0$ such that every $\mathcal{F}$-free graph $G$ on $n\geq N$ vertices with $\delta(G)\geq (\pi(\mathcal{F})-\varepsilon)n$ is a subgraph of some member of $\mathcal{G}$, then we say that $\mathcal{F}$ is \textit{degree-stable} with respect to $\mathcal{G}$.

Given a family of graphs   $\mathcal{G}$, let  $\mathcal{G}_{n}$  denote the collection of all    $n$-vertex graphs contained in   $\mathcal{G}$. For   $0\leq \alpha <1$, we define
$$ \lambda_{\alpha}(\mathcal{G}_n)=\max_{G\in \mathcal{G}_{n}} \lambda_{\alpha}(G).$$

Nikiforov \cite{N1} established the following results, which will be invoked in our proofs.
\begin{lem} [Nikiforov \cite{N1}] \label{c1}
    Let $G$ be a graph and $0\leq \alpha\leq 1$. Then
     $$\alpha \Delta(G)\leq \lambda_{\alpha}(G) \leq \alpha \Delta(G)+(1-\alpha)\lambda (G).$$
\end{lem}

\begin{lem}[Nikiforov \cite{N1}] \label{LL5}
    Let $0\leq \alpha \leq 1$. If $G$ is of order $n$ and $m$ edges, then
  $$\lambda_{\alpha} (G)\geq \sqrt{\frac{1}{n}\sum_{u\in V(G)}d^2(u)} \ \text { and }\  \lambda_{\alpha}(G)\geq \frac{2m}{n} .$$
  Equality holds in the second inequality if and only if $G$ is regular. If $\alpha>0$, equality holds in the first inequality if and only if $G$ is regular.
\end{lem}

\begin{lem}\label{T1}
 For $r\geq 2$ and $0\leq \alpha \leq 1-\frac{1}{r}$, let $G$ be a graph of order $n$ and let $\mathbf{x}=(x_{1}, \dots, x_{n})$ be a non-negative unit eigenvector of $A_{\alpha}(G)$ corresponding to $\lambda_{\alpha}(G)$. Suppose $w\in V(G)$  satisfies $x_{w}=\min \{x_{1}, \dots, x_{n}\}$. Then
    $$\lambda_{\alpha}(G-w)\geq \lambda_{\alpha}(G)\frac{1-2x^2_{w}}{1-x^2_{w}}-\alpha\frac{1-nx^2_{w}}{1-x^2_{w}}.$$

\end{lem}

\noindent\textbf{Proof.}
For simplicity, we write $x=x_w$. By the  Rayleigh principle, it follows that
    \begin{eqnarray}\label{1}
		\lambda_{\alpha}(G)&=& \mathbf{x}^{T}A_{\alpha}(G)\mathbf{x} \nonumber\\[2mm]
           &=& \sum_{uv\in E(G)}\left(\alpha x^2_{u}+2(1-\alpha) x_{u}x_{v}+\alpha x^2_{v} \right)\nonumber \\[2mm]
           &=& \sum_{uv\in E(G-w)}\left(\alpha x^2_{u}+2(1-\alpha) x_{u}x_{v}+\alpha x^2_{v} \right)+\sum_{u\in N(w)} \left(\alpha x^2+2(1-\alpha)xx_u+\alpha x^2_{u} \right)\nonumber \\[2mm]
           &\leq & \left(1-x^2\right)\lambda_{\alpha}(G-w)+\alpha d(w)x^2+ 2(1-\alpha)x \sum_{u\in N(w)}x_{u}+\alpha \sum_{u\in N(w)} x^2_{u}.		
    \end{eqnarray}

Note that $(1-\alpha)\sum_{u\in N(w)}x_{u}=(\lambda_{\alpha}(G)-\alpha d(w))x$. Consequently, we obtain the following inequality:
 \begin{eqnarray*}
	&	&\alpha d(w)x^2+ 2(1-\alpha)x \sum_{u\in N(w)}x_{u}+\alpha \sum_{u\in N(w)} x^2_{u}\\[2mm]
             &\leq& \alpha d(w)x^2+2(\lambda_{\alpha}(G)-\alpha d(w))x^2+\alpha(1-(n-d(w))x^2)\\[2mm]
              &=&2\lambda_{\alpha}(G)x^2-\alpha nx^2+\alpha.
    \end{eqnarray*}

Combining this with (\ref{1}), we obtain
$$\lambda_{\alpha}(G-w)\geq \lambda_{\alpha}(G)\frac{1-2x^2}{1-x^2}-\alpha\frac{1-nx^2}{1-x^2}.$$
\qed

\begin{lem} \label{T2}
    Let $\mathcal{G}$ be a hereditary graph family. For $r\geq 2$ and $0\leq \alpha \leq 1-\frac{1}{r}$ , if $\lambda_{\alpha}(\mathcal{G}_n)> \left(1-\frac{1}{r}\right)n-\left(1-\frac{1}{r}\right)$, then the limit
    $$\pi_{\alpha}(\mathcal{G}):=\lim_{n\to \infty} \frac{\lambda_{\alpha}(\mathcal{G}_n)}{n}$$
    exists, and it satisfies
    $$\pi_{\alpha}(\mathcal{G})\leq \frac{\lambda_{\alpha}(\mathcal{G}_n)}{n-1}.$$
\end{lem}

\noindent\textbf{Proof.}
Let $G\in \mathcal{G}_n$ satisfy $\lambda_{\alpha}(G)=\lambda_{\alpha}(\mathcal{G}_n)$, and let $\mathbf{x}=(x_{1}, \dots, x_{n})$ be a non-negative unit eigenvector of $A_{\alpha}(G)$ corresponding to $\lambda_{\alpha}(G)$. Suppose $w\in V(G)$ such  that $x_{w}=\min \{x_{1}, \dots, x_{n}\}$.
Since  $\mathcal{G}$  is hereditary, $G-w\in \mathcal{G}_{n-1}$,
which implies that $\lambda_{\alpha}(\mathcal{G}_{n-1})\geq \lambda_{\alpha}(G-w)$. For simplicity, let $x=x_w$. Note that $x^2\leq \frac{1}{n}$. By Lemma \ref{T1}, we obtain

$$\lambda_{\alpha}(\mathcal{G}_{n-1})\geq \lambda_{\alpha}(\mathcal{G}_n)\frac{1-2x^2}{1-x^2}-\alpha\frac{1-nx^2}{1-x^2}.$$
Thus,
\begin{eqnarray*}
	\frac{\lambda_{\alpha}(\mathcal{G}_{n-1})}{n-2}
    &\geq& \frac{\lambda_{\alpha}(\mathcal{G}_{n})}{n-1}\left(1+\frac{1}{n-2}\right)\frac{1-2x^2}{1-x^2}-\alpha\frac{1-nx^2}{(n-2)(1-x^2)} \nonumber\\[2mm]
    &=& \frac{\lambda_{\alpha}(\mathcal{G}_{n})}{n-1}\left(1+\frac{1-nx^2}{(n-2)(1-x^2)}\right)-\alpha\frac{1-nx^2}{(n-2)(1-x^2)} \nonumber\\[2mm]
    &\geq & \frac{\lambda_{\alpha}(\mathcal{G}_{n})}{n-1},
    \end{eqnarray*}
    where the last inequality holds as $\lambda_{\alpha}(\mathcal{G}_n)> \left(1-\frac{1}{r}\right)n-\left(1-\frac{1}{r}\right)$ and $0\leq \alpha \leq 1-\frac{1}{r}$.

Therefore, $\frac{\lambda_{\alpha}(\mathcal{G}_{n-1})}{n-2}\geq \frac{\lambda_{\alpha}(\mathcal{G}_{n})}{n-1}>0$. This implies that  the limit  $\pi_{\alpha}(\mathcal{G}):=\lim_{n\to \infty} \frac{\lambda_{\alpha}(\mathcal{G}_n)}{n}$ exists.  Moreover, $\pi_{\alpha}(\mathcal{G})\leq \frac{\lambda_{\alpha}(\mathcal{G}_n)}{n-1}$. \qed

\section {Proofs of Theorems \ref{TT1} and \ref{TT11}}

\begin{lem} \label{T4}
   For $r\geq 2$ and $0\leq \alpha\leq 1-\frac{1}{r}$, let $G$ be a graph on $n$ vertices, and  $\mathbf{x}=(x_{1}, \dots, x_{n})$ be a non-negative unit eigenvector of $A_{\alpha}(G)$ corresponding to $\lambda_{\alpha}(G)$. Denote $x=\min \{x_{1}, \dots, x_{n}\}$ and $\delta=\delta(G)$. If $x>0$, then
    $$\lambda_{\alpha}(G)\leq \alpha \delta+ (1-\alpha)\sqrt{\delta^2+\left(\frac{1}{nx^2}-1\right)n\delta}.$$
\end{lem}

\noindent\textbf{Proof.}
Let $u\in V(G)$ such that $d(u)=\delta$. According to the eigenvector equation with respect to $u$, we obtain $$\frac{1}{1-\alpha}(\lambda_{\alpha}(G)-\alpha \delta)x_{u}=\sum_{v\in N(u)}x_{v}.$$ Hence,
\begin{eqnarray*}
\left(\frac{1}{1-\alpha}\right)^2 \left(\lambda_{\alpha}(G)-\alpha \delta\right)^2x^2&\leq&\left(\frac{1}{1-\alpha}\right)^2 \left(\lambda_{\alpha}(G)-\alpha \delta\right)^2x^2_{u}\\[2mm]
&\leq& \delta\sum_{v\in N(u)}x^2_{v}\\[2mm]
&\leq&  \delta(1-(n-\delta)x^2),
\end{eqnarray*}
which implies that $$\lambda_{\alpha}(G)\leq \alpha \delta+ (1-\alpha)\sqrt{\delta^2+\left(\frac{1}{nx^2}-1\right)n\delta}.$$ \qed

Let $\mathcal{F}$ be a family of graphs with $\chi(\mathcal{F})=r+1\geq 3$, and let $\mathcal{G}$ denote the collection of all $\mathcal{F}$-free graphs. Clearly, $\mathcal{G}$ is a hereditary graph family,  and $\lambda_{\alpha}(\mathcal{G}_n)=\mathrm{ex}_{\alpha}(n,\mathcal{F})$. Note that the   $r$-partite Tur\'an graph
$T_{n,r}\in \mathcal{G}_n$. By Lemma \ref{LL5}, we thus have
$$\lambda_{\alpha}(\mathcal{G}_n)\geq \lambda_{\alpha}(T_{n,r})\geq \frac{2e(T_{n,r})}{n}\geq \left(1-\frac{1}{r}\right)n-\frac{r}{4n}.$$
 Furthermore, Lemma \ref{T2} guarantees the existence of the limit
  $\lim_{n \to \infty } \frac{\mathrm{ex}_{\alpha}(n,\mathcal{F})}{n}$. We denote this limit as $\pi'_{\alpha}(\mathcal{F})$, so that $\pi'_{\alpha}(\mathcal{F})=\pi_{\alpha}(\mathcal{G})$.

\begin{lem} \label{T5}
  For $r\geq 2$ and $0\leq \alpha\leq 1- \frac{1}{r}$, let $\mathcal{F}$ be a family of graphs with $\chi(\mathcal{F})=r+1\geq 3$.  Then $\pi'_{\alpha}(\mathcal{F})\geq \pi(\mathcal{F})$. Furthermore, if $\pi'_{\alpha}(\mathcal{F})>1-\frac{1}{r}$, then $\pi'_{\alpha}(\mathcal{F})=\pi(\mathcal{F})$.
\end{lem}

\noindent\textbf{Proof.}
We first prove that $\pi'_{\alpha}(\mathcal{F})\geq \pi(\mathcal{F})$. Let $G\in \mathrm{Ex}(n,\mathcal{F})$. By Lemma \ref{LL5},
$$\frac{\mathrm{ex}_{\alpha}(n,\mathcal{F})}{n}\geq \frac{\lambda_{\alpha}(G)}{n}\geq \frac{2\mathrm{ex}(n,\mathcal{F})}{n^2},$$
which implies that $\pi'_{\alpha}(\mathcal{F})\geq \pi(\mathcal{F})$.

Now, suppose $\pi'_{\alpha}(\mathcal{F})>1-\frac{1}{r}$. We will show $\pi'_{\alpha}(\mathcal{F})\leq \pi(\mathcal{F})$, thus finishing the proof.
Clearly, $\pi'_{\alpha}(\mathcal{F})<1$. So we assume $\pi'_{\alpha}(\mathcal{F})=1-\frac{1}{r}+\varepsilon$, where $0<\varepsilon<\frac{1}{r}$.
Let $\mathcal{G}$ be the collection of all $\mathcal{F}$-free graphs. Let $G_{n}\in \mathcal{G}_{n}$ such that $\lambda_{\alpha}(G_n)=\lambda_{\alpha}(\mathcal{G}_{n})$. Suppose $\mathbf{x}^{(n)}=\left(x^{(n)}_{1}, \dots, x^{(n)}_{n}\right)$ is a non-negative unit eigenvector of $A_{\alpha}(G_n)$ corresponding to $\lambda_{\alpha}(G_n)$. Denote $x^{(n)}=\min \left\{x^{(n)}_{1}, \dots, x^{(n)}_{n} \right\}$.

\begin{claim}
  There exist infinitely many $n\in \mathbb{N}$ such that
  $$ \left(x^{(n)}\right)^2>\frac{1}{n}\left(1-\frac{1}{\ln{n}}\right)>0. $$
\end{claim}

\noindent\textbf{Proof.}
Suppose, for contradiction, that there exists some  $n_0$ such that for all $n\geq n_0$,
$$0\leq \left(x^{(n)}\right)^2\leq \frac{1}{n}\left(1-\frac{1}{\ln{n}}\right).$$
Note that $\mathcal{G}$ is a hereditary graph family. By Lemma \ref{T1}, we thus have
    $$\lambda_{\alpha}(\mathcal{G}_{n-1})\geq \lambda_{\alpha}(\mathcal{G}_n)\frac{1-2\left(x^{(n)}\right)^2}{1-\left(x^{(n)}\right)^2}-\alpha\frac{1-n\left(x^{(n)}\right)^2}{1-\left(x^{(n)}\right)^2}.$$
Moreover, since $\lambda_{\alpha}(\mathcal{G}_n)\geq \lambda_{\alpha}(T_{n,r})\geq \left(1-\frac{1}{r}\right)n-\frac{r}{4n}$, it follows from Lemma \ref{T2} that $\frac{\lambda_{\alpha}(\mathcal{G}_n)}{n-1}\geq \pi_{\alpha}(\mathcal{G})=\pi'_{\alpha}(\mathcal{F})=1-\frac{1}{r}+\varepsilon$. Now, let $n\geq n_0$ and let $x=x^{(n)}$.
Then we get
\begin{eqnarray*}
	\frac{\lambda_{\alpha}(\mathcal{G}_{n-1})}{n-2}
    &\geq& \frac{\lambda_{\alpha}(\mathcal{G}_{n})}{n-1}\left(1+\frac{1}{n-2}\right)\frac{1-2x^2}{1-x^2}-\alpha\frac{1-nx^2}{(n-2)(1-x^2)} \nonumber\\[2mm]
    &=& \frac{\lambda_{\alpha}(\mathcal{G}_{n})}{n-1}\left(1+\left(1-\alpha\frac{n-1}{\lambda_{\alpha}(\mathcal{G}_{n})}\right)\frac{1-nx^2}{(n-2)(1-x^2)}\right) \nonumber\\[2mm]
    &\geq & \frac{\lambda_{\alpha}(\mathcal{G}_{n})}{n-1}\left(1+ \frac{\varepsilon(1-nx^2)}{(n-2)(1-x^2)}\right)\nonumber\\[2mm]
     &\geq & \frac{\lambda_{\alpha}(\mathcal{G}_{n})}{n-1}\left(1+ \frac{\varepsilon}{n \ln n}\right)\nonumber,
\end{eqnarray*}
     where the second inequality holds as $0\leq \alpha\leq 1- \frac{1}{r}$ and $\frac{n-1}{\lambda_{\alpha}(\mathcal{G}_n)}\leq \frac{1}{1-\frac{1}{r}+\varepsilon}$,  and the last inequality holds as $x^2\leq \frac{1}{n}\left(1-\frac{1}{\ln{n}}\right)$. Thus, for any $i\geq n_0$, the following inequality holds:
     $$ \frac{\lambda_{\alpha}(\mathcal{G}_{i-1})}{i-2}-  \frac{\lambda_{\alpha}(\mathcal{G}_{i})}{i-1}\geq \frac{\varepsilon}{i \ln i}\cdot\frac{\lambda_{\alpha}(\mathcal{G}_{i})}{i-1}.$$
     Summing both sides over  $i$  from  $n_0$  to  $2n$, we further obtain
\begin{eqnarray*}
 \frac{\lambda_{\alpha}(\mathcal{G}_{n_0-1})}{n_0-2}-  \frac{\lambda_{\alpha}(\mathcal{G}_{2n})}{2n-1}
 &=& \sum_{i=n_0}^{2n} \left(\frac{\lambda_{\alpha}(\mathcal{G}_{i-1})}{i-2}-  \frac{\lambda_{\alpha}(\mathcal{G}_{i})}{i-1} \right) \nonumber\\[2mm]
 &\geq& \sum_{i=n_0}^{2n} \frac{\varepsilon}{i \ln i}\cdot\frac{\lambda_{\alpha}(\mathcal{G}_{i})}{i-1} \nonumber\\[2mm]
 &\geq& \varepsilon\pi'_{\alpha}(\mathcal{F}) \sum_{i=n_0}^{2n}\frac{1}{i \ln i},
\end{eqnarray*}
where the last inequality holds as $\frac{\lambda_{\alpha}(\mathcal{G}_i)}{i-1}\geq \pi_{\alpha}(\mathcal{G})=\pi'_{\alpha}(\mathcal{F})$. Taking $n$ sufficiently large, we get a contradiction as $\sum_{i=n_0}^{2n}\frac{1}{i \ln i}$ diverges. This completes the proof of Claim 1. \qed

By Claim 1, there exists an infinite sequence $\{n_k\}_{k=1}^\infty$ of positive integers such that for each $n\in \{n_k\}_{k=1}^\infty$,
$$\left(x^{(n)}\right)^2>\frac{1}{n}\left(1-\frac{1}{\ln{n}}\right)>0.$$
In light of Lemma \ref{T4}, we deduce that
\begin{eqnarray*}
\frac{\lambda_{\alpha}(G_n)}{n}
&\leq& \frac{\alpha \delta(G_n)}{n}+ (1-\alpha)\sqrt{\left(\frac{\delta(G_n)}{n}\right)^2+\left(\frac{1}{n \left(x^{(n)}\right)^2}-1\right)\frac{\delta(G_n)}{n}}\nonumber\\[2mm]
&\leq& \frac{\alpha \delta(G_n)}{n} + (1-\alpha)\left(\frac{\delta(G_n)}{n}+\frac{1}{1-\frac{1}{\ln n}}-1\right) \nonumber\\[2mm]
&\leq & \frac{2\mathrm{ex}(n,\mathcal{F})}{n^2}+(1-\alpha)\left(\frac{1}{1-\frac{1}{\ln n}}-1\right).
\end{eqnarray*}
Therefore, we obtain
$$\pi'_{\alpha}(\mathcal{F})= \lim_{n\to \infty}\frac{\lambda_{\alpha}(G_n)}{n}\leq \lim_{n\to \infty}\left(\frac{2\mathrm{ex}(n,\mathcal{F})}{n^2}+(1-\alpha)\left(\frac{1}{1-\frac{1}{\ln n}}-1\right)\right)=\pi(\mathcal{F}).$$
Thus, $\pi'_{\alpha}(\mathcal{F})=\pi(\mathcal{F})$. The proof is complete. \qed

We are now ready to prove Theorem \ref{TT1}.

\noindent\textbf{Proof of Theorem \ref{TT1}.}
 It suffices to show that $\pi'_{\alpha}(\mathcal{F})=1-\frac{1}{r}$. We prove the assertion by contradiction. Recall that $\pi'_{\alpha}(\mathcal{F})\geq\pi(\mathcal{F})=1-\frac{1}{r}$. If $\pi'_{\alpha}(\mathcal{F})> 1-\frac{1}{r}$, then by Lemma \ref{T5}, we deduce that $\pi'_{\alpha}(\mathcal{F})=\pi(\mathcal{F})=1-\frac{1}{r}$, which  is a contradiction. Thus $\pi'_{\alpha}(\mathcal{F})=\pi(\mathcal{F})=1-\frac{1}{r}$. This completes the proof. \qed

Next, we give the proof of Theorem \ref{TT11}.

\noindent\textbf{Proof of Theorem \ref{TT11}.}
 We first consider the case where $\mathcal{F}$   does not contain  $K_{1,k}\cup s K_1$  for any   $k\geq 1, s\geq 0$. In this case, $K_{1,n-1}$   is   $\mathcal{F}$-free. Let   $G\in \mathrm{Ex}_{\alpha}(n,\mathcal{F})$. By Lemma \ref{c1}, we have   $\lambda_{\alpha}(G)\geq \lambda_{\alpha}(K_{1,n-1})\geq\alpha (n-1)$. Since   $\chi(\mathcal{F})=2$, there exist    a number $t>0$ and $F\in \mathcal{F}$ such that $F\subset K_{t,t}$. By the well-known K\H{o}v\'ari-S\'{o}s-Tur\'{a}n
  theorem,   $\mathrm{ex}(n,K_{t,t})=O(n^{2-\frac{1}{t}})$.
 Therefore, by Lemma \ref{c1}, we obtain
 $$\alpha (n-1)\leq \lambda_{\alpha}(G) \leq \alpha \Delta(G)+(1-\alpha)\lambda(G)\leq  \alpha (n-1)+ (1-\alpha)\sqrt{2e(G)}=  \alpha n+ o(n).$$
Then  $\mathrm{ex}_{\alpha}(n,\mathcal{F})=(\alpha +o(1))n$. Next, we assume $\mathcal{F}$ contains $K_{1,k}\cup s K_1$ for some $k\geq 1, s\geq 0$. Then $\Delta(G)\leq k-1$. It follows from Lemma \ref{c1} that
$$\lambda_{\alpha}(G) \leq \alpha \Delta(G)+(1-\alpha)\lambda(G)\leq \Delta(G)\leq k-1.$$
Thus, $\mathrm{ex}_{\alpha}(n,\mathcal{F})=O(1)$. If $\mathcal{F}=\{K_{1,k}\}$, then  $K_{k} \cup I_{n-k}$ is $\mathcal{F}$-free. Thus $\lambda_{\alpha}(G) \geq \lambda_{\alpha}(K_{k} \cup I_{n-k})=k-1$. Hence, $\mathrm{ex}_{\alpha}(n,\mathcal{F})=k-1$. The proof  is complete.
\qed

\section{Proof of Theorem \ref{TT2} }

\subsection{Some auxiliary lemmas}

Keevash, Lenz, and Mubayi \cite{PK} established a criterion that allows spectral extremal problems involving the   $p$-spectral radius of hypergraphs to be deduced from their corresponding hypergraph Turán problems. Zheng, Li, and Li \cite{Z2} later proposed a similar criterion for the signless Laplacian spectral radius. Motivated by these works, we develop a criterion for the   $\alpha$-spectral radius, which reduces   $\alpha$-spectral Tur\'an-type problems to classical extremal problems satisfying the degree-stable property.

\begin{thm}\label{TT3}
    Let $\mathcal{F}$ be a family of graphs with $\chi(\mathcal{F})=r+1\geq 3$. Suppose $0< \varepsilon < \frac{1}{2}$, $\sigma<\frac{\varepsilon^3}{7}$ and $0\leq \alpha\leq 1-\frac{1}{r}-\varepsilon$. Let $\mathcal{G}'_n$ be the set of all $\mathcal{F}$-free graphs on $n$ vertices with minimum degree larger than $(\pi(\mathcal{F})-\varepsilon)n$. Suppose that there exists $N>0$ such that for every $n>N$, we have
\begin{eqnarray} \label{2}
 |\mathrm{ex}(n,\mathcal{F})-\mathrm{ex}(n-1,\mathcal{F})-\pi(\mathcal{F})n|\leq \sigma n
\end{eqnarray}
and
\begin{eqnarray} \label{3}
 \left|\lambda_{\alpha}(\mathcal{G}'_{n})-\frac{2\mathrm{ex}(n,\mathcal{F})}{n}\right|\leq \sigma.
\end{eqnarray}
Then there exists $n_{0}>0$ such that for any $\mathcal{F}$-free graph $G$ on $n\geq n_0$ vertices, we have
$$\lambda_{\alpha}(G)\leq \lambda_{\alpha}(\mathcal{G}'_{n}).$$
In addition, if the equality holds, then $G\in \mathcal{G}'_{n}$.

\end{thm}

Now we  prove  Theorem \ref{TT3}. Let $G$ be an $\mathcal{F}$-free graph on $n$ vertices. Let $\mathbf{x}=(x_{1}, \dots, x_{n})$ be a non-negative unit eigenvector of $A_{\alpha}(G)$ corresponding to $\lambda_{\alpha}(G)$. Suppose $w\in V(G)$ such that $x=x_{w}=\min \{x_{1}, \dots, x_{n}\}$. Let $\delta=\delta(G)$.

From \eqref{3}, it is readily seen that   $\frac{ \lambda_{\alpha}(\mathcal{G}'_{n})}{n}-\frac{2\mathrm{ex}(n,\mathcal{F})}{n^2}=o(1)$. Recalling that  $\pi(\mathcal{F})=\lim_{n\to \infty}\frac{2\mathrm{ex}(n,\mathcal{F})}{n^2}$, this immediately implies
\begin{equation}\label{0}
\lambda_{\alpha}(\mathcal{G}'_{n})=(\pi(\mathcal{F})+o(1))n.
\end{equation}
Clearly,   $\delta(T_{n,r})\geq \left(1-1/r\right)n-1>\left(1-1/r-\varepsilon\right)n$, and   $T_{n,r}$ is $\mathcal{F}$-free. Thus,   $T_{n,r}\in \mathcal{G}'_n$. Combining Lemma \ref{LL5} with the known lower bound   $e(T_{n,r})\geq \frac{r-1}{2r}n^2-\frac{r}{8}$,
we deduce
\begin{eqnarray}\label{4}
    \lambda_{\alpha}( \mathcal{G}'_n)\geq \left(1-\frac{1}{r}\right)n-\frac{r}{4n}.
\end{eqnarray}

\begin{lem} \label{L1}
    $\lambda_{\alpha}(\mathcal{G}'_n)\geq \lambda_{\alpha}( \mathcal{G}'_{n-1})+\pi(\mathcal{F})-5\sigma$.
\end{lem}

\noindent\textbf{Proof.}
By the triangle inequality, combining \eqref{2} and \eqref{3}, we deduce that
\begin{eqnarray*}
   & &\left | \lambda_{\alpha}( \mathcal{G}'_n)- \lambda_{\alpha}( \mathcal{G}'_{n-1})-\pi(\mathcal{F})\right|\\[2mm]
    &\leq& \left|\frac{2\mathrm{ex}(n,\mathcal{F})}{n}-\frac{2\mathrm{ex}(n-1,\mathcal{F})}{n-1}-\pi(\mathcal{F}) \right| +2\sigma \nonumber \\[2mm]
   &=& \left|\frac{2}{n}(\mathrm{ex}(n,\mathcal{F})-\mathrm{ex}(n-1,\mathcal{F})-\pi (\mathcal{F})n)+ \pi (\mathcal{F})-\frac{2\mathrm{ex}(n-1,\mathcal{F})}{n(n-1)}\right| +2\sigma \nonumber \\[2mm]
   &\leq & \left|\pi(\mathcal{F})-\frac{2\mathrm{ex}(n-1,\mathcal{F})}{n(n-1)} \right|+4\sigma \nonumber.
\end{eqnarray*}
Note that $\lim_{n\to \infty}\frac{2\mathrm{ex}(n-1,\mathcal{F})}{n(n-1)}=\pi(\mathcal{F})$. Thus, for sufficiently large $n$, we have
\begin{equation*}
    \begin{aligned}
 \lambda_{\alpha}( \mathcal{G}'_n)\geq \lambda_{\alpha}( \mathcal{G}'_{n-1})+\pi(\mathcal{F})-5\sigma.
    \end{aligned}
\end{equation*}
\qed

\begin{lem} \label{L2}
 $x^2 \leq \frac{\delta (1-\alpha)^2}{(\lambda_{\alpha}(G)-\alpha \delta)^2+\delta(n-\delta) (1-\alpha)^2}$.
\end{lem}

\noindent\textbf{Proof.}
Let $u\in V(G)$ such that $d(u)=\delta$. According to the eigenvector equation with respect to the vertex $u$, we get
$$ (\lambda_{\alpha}(G)-\alpha \delta)x_{u}=(1-\alpha)\sum_{v\in N(u)}x_{v}.$$
Combining this with Power Mean inequality, we obtain
\begin{eqnarray*}
\left(\frac{1}{1-\alpha}(\lambda_{\alpha}(G)-\alpha \delta)x\right)^2&\leq& \left(\sum_{v\in N(u)}x_{v}\right)^2
\leq \delta \sum_{v\in N(u)}x^2_{v}\\[2mm]
&\leq& \delta\left(1-(n-\delta)x^2\right),
\end{eqnarray*}
which implies that
$$x^2 \leq \frac{\delta (1-\alpha)^2}{(\lambda_{\alpha}(G)-\alpha \delta)^2+\delta(n-\delta) (1-\alpha)^2}.$$ \qed

\begin{lem} \label{L3}
 If $\lambda_{\alpha}(G)\geq \lambda_{\alpha}(\mathcal{G}'_{n})$ and  $\delta\leq (\pi (\mathcal{F})-\varepsilon)n$, then for sufficiently large $n$, we have
 $$x^2<\frac{1-\varepsilon^2}{n}.$$
\end{lem}

\noindent\textbf{Proof.}
 If $\delta=0$, then  $x=0$, the result follows. In the following, we assume $\delta\geq 1$.
 Set $\theta=\frac{\varepsilon}{2}$. By the assumption that $\lambda_{\alpha}(G)\geq \lambda_{\alpha}(\mathcal{G}'_{n})$ and (\ref{4}), we get
 \begin{eqnarray}\label{eq6}
 \lambda_{\alpha}(G)\geq \lambda_{\alpha}( \mathcal{G}'_n)\geq \left(1-\frac{1}{r}\right)n-\frac{r}{4n}\geq (\pi(\mathcal{F})-\theta)n.
  \end{eqnarray}
   Let $$f(y,z)=\frac{1}{y}(z-\alpha y)^2+(n-y) (1-\alpha)^2,$$  where $1\leq y\leq (\pi(\mathcal{F})-\varepsilon)n$ and $z\geq (\pi(\mathcal{F})-\theta)n$.

 A straightforward calculation shows
  $$\frac{\partial f(y,z)}{\partial y}=-\left(\frac{z^2}{y^2}+1-2\alpha\right)  \mbox{and} \ \frac{\partial f(y,z)}{\partial z}=2\left(\frac{z}{y}-\alpha\right),$$ which implies that $\frac{\partial f(y,z)}{\partial y}\leq 0$ and $\frac{\partial f(y,z)}{\partial z}> 0$ for $1\leq y\leq (\pi(\mathcal{F})-\varepsilon)n$ and $z\geq (\pi(\mathcal{F})-\theta)n$.

 Combining this with (\ref{eq6}) and assumption $\delta\leq (\pi (\mathcal{F})-\varepsilon)n$, we get
 $$f(\delta,\lambda_{\alpha}(G))\geq f((\pi(\mathcal{F})-\varepsilon)n,(\pi(\mathcal{F})-\theta)n).$$

 By Lemma \ref{L2}, we have
 \begin{eqnarray*}
   x^2n &\leq& \frac{(1-\alpha)^2}{f(\delta,\lambda_{\alpha}(G))}n\leq \frac{(1-\alpha)^2}{f((\pi(\mathcal{F})-\varepsilon)n,(\pi(\mathcal{F})-\theta)n)}n\nonumber\\[2mm]
    &=& \frac{(1-\alpha)^2(\pi(\mathcal{F})-\varepsilon)}{((1-\alpha)(\pi(\mathcal{F})-\varepsilon)+\varepsilon-\theta)^2+(1-\alpha)^2(\pi(\mathcal{F})-\varepsilon)(1-\pi(\mathcal{F})+\varepsilon)} \nonumber\\[2mm]
   &=& \frac{(1-\alpha)^2(\pi(\mathcal{F})-\varepsilon)}{(1-\alpha)^2(\pi(\mathcal{F})-\varepsilon)+2(1-\alpha)(\pi(\mathcal{F})-\varepsilon)(\varepsilon-\theta)+(\varepsilon-\theta)^2} \\[2mm]
    &\leq & \frac{1}{1+\frac{2(\varepsilon-\theta)}{1-\alpha}}< 1- \varepsilon^2,
\end{eqnarray*}
completing the proof.
 \qed

\begin{lem}\label{L4}
    Suppose $G$ satisfies $\lambda_{\alpha}(G)\geq \lambda_{\alpha}(\mathcal{G}'_{n})$ and $x^2<\frac{1-\varepsilon^2}{n}$. Then, for sufficiently large $n$, we have
\begin{equation}\label{5}
   \lambda_{\alpha}(G-w)\geq \lambda_{\alpha}(G)\left(1-\frac{1-\varepsilon^3}{n-1}\right)
\end{equation}
     and
\begin{equation}\label{6}
     \lambda_{\alpha}(G-w)>\lambda_{\alpha}(\mathcal{G}'_{n-1}).
\end{equation}

\end{lem}

\noindent\textbf{Proof.}
By Lemma \ref{T1}, we obtain
\begin{equation*}
    \begin{aligned}
 \lambda_{\alpha}(G-w)\geq \lambda_{\alpha}(G)\frac{1-2x^2}{1-x^2}-\alpha\frac{1-nx^2}{1-x^2}.
    \end{aligned}
\end{equation*}
Combining this with the assumption $\lambda_{\alpha}(G)\geq \lambda_{\alpha}(\mathcal{G}'_{n})$, (\ref{4}), and $\alpha\leq 1-\frac{1}{r}-\varepsilon$, we get
\begin{eqnarray}
\frac{\lambda_{\alpha}(G-w)}{n-2}
    &\geq& \frac{\lambda_{\alpha}(G)}{n-1}\left(1+\frac{1}{n-2}\right)\frac{1-2x^2}{1-x^2}-\alpha\frac{1-nx^2}{(n-2)(1-x^2)} \label{7}\\[2mm]
    &=& \frac{\lambda_{\alpha}(G)}{n-1}\left(1+\left(1-\alpha\frac{n-1}{\lambda_{\alpha}(G)}\right)\frac{1-nx^2}{(n-2)(1-x^2)}\right) \nonumber\\[2mm]
    &\geq & \frac{\lambda_{\alpha}(G)}{n-1}\left(1+\left(1-\left(1-\frac{1}{r}-\varepsilon\right)\frac{n-1}{\left(1-\frac{1}{r}\right)n-\frac{r}{4n}}\right)\cdot \frac{1-nx^2}{(n-2)(1-x^2)}\right)\nonumber\\[2mm]
    &\geq & \frac{\lambda_{\alpha}(G)}{n-1}\left(1+\frac{\varepsilon}{1-\frac{1}{r}}\cdot \frac{1-nx^2}{(n-2)(1-x^2)}\right)\nonumber\\[2mm]
    &\geq & \frac{\lambda_{\alpha}(G)}{n-1} \left(1+ \frac{r\varepsilon^3}{(r-1)(n-1)}\right)\nonumber,
\end{eqnarray}
where the last inequality holds as $x^2<\frac{1-\varepsilon^2}{n}$.

Then,  it follows that
\begin{eqnarray*}
\lambda_{\alpha}(G-w)&\geq& \lambda_{\alpha}(G)\left(1-\frac{1}{n-1}\right)\left(1+ \frac{r\varepsilon^3}{(r-1)(n-1)}\right)\\[2mm]
&\geq & \lambda_{\alpha}(G)\left(1-\frac{1-\varepsilon^3}{n-1}\right).
\end{eqnarray*}

Next we prove (\ref{6}). By (\ref{7}), we obtain
\begin{eqnarray}
    \frac{\lambda_{\alpha}(G-w)}{n-2}
     &\geq & \frac{\lambda_{\alpha}(G)}{n-1}\left(1+\frac{1-nx^2}{(n-2)(1-x^2)}\right)-\alpha\frac{1-nx^2}{(n-2)(1-x^2)}\nonumber\\[2mm]
     &=& \frac{\lambda_{\alpha}(G)}{n-1} + \left(\frac{\lambda_{\alpha}(G)}{n-1}-\alpha \right)\frac{1-nx^2}{(n-2)(1-x^2)}\nonumber\\[2mm]
     &\geq & \frac{\lambda_{\alpha}(G)}{n-1} + \left(\frac{\left(1-\frac{1}{r}\right)(n-1)}{n-1}-\left(1-\frac{1}{r}-\varepsilon\right) \right)\frac{1-nx^2}{(n-2)(1-x^2)}\nonumber\\[2mm]
     &= &  \frac{\lambda_{\alpha}(G)}{n-1} +\frac{\varepsilon(1-nx^2)}{(n-2)(1-x^2)}\nonumber\\[2mm]
     &> & \frac{\lambda_{\alpha}(G)}{n-1} + \frac{\varepsilon^3}{n-2},\label{eq10}
\end{eqnarray}
where the second inequality holds as $\lambda_{\alpha}(G)\geq \left(1-\frac{1}{r}\right)(n-1)$ and $\alpha\leq 1-\frac{1}{r}-\varepsilon$, and the last inequality holds as $x^2<\frac{1-\varepsilon^2}{n}$.

For sufficiently large $n$, by (\ref{0}), we have $ \lambda_{\alpha}(\mathcal{G}'_{n-1})\leq (\pi(\mathcal{F})+\sigma)(n-1).$
Combining this with (\ref{eq10}), Lemma \ref{L1} and the assumption $\lambda_{\alpha}(G)\geq \lambda_{\alpha}(\mathcal{G}'_{n})$, we get
\begin{eqnarray*}
    \lambda_{\alpha}(G-w)&\geq& \lambda_{\alpha}(G)\left(1-\frac{1}{n-1}\right)+\varepsilon^3\\
     &\geq &(\lambda_{\alpha}( \mathcal{G}'_{n-1})+\pi(\mathcal{F})-5\sigma)\left(1-\frac{1}{n-1}\right)+\varepsilon^3\\[2mm]
     &\geq & \lambda_{\alpha}( \mathcal{G}'_{n-1})-7\sigma+\varepsilon^3 \\[2mm]
     &>&\lambda_{\alpha}( \mathcal{G}'_{n-1}),
\end{eqnarray*}
where the last inequality holds as $\sigma<\frac{\varepsilon^3}{7}$. The proof is complete. \qed

 The following facts, proved by Zheng, Li,  and Li \cite{Z2}, will be employed in our proof

 \vspace{2mm}
\noindent\textbf{Fact 1.} If $0<x<\frac{1}{2}$ and $0<a<1$, then $\ln (1-ax)+ax+x^2>0$.

\vspace{2mm}
\noindent\textbf{Fact 2.} If $x>1$, then $\frac{1}{x}<\ln x-\ln(x-1)$ and $\frac{1}{x^2}<\frac{1}{x-1}-\frac{1}{x}$.
\vspace{2mm}


We are now ready to prove Theorem \ref{TT3}.

\noindent\textbf{Proof of Theorem \ref{TT3}.}
We prove by contradiction. Suppose to the contrary that there exists an $\mathcal{F}$-free graph $G$ with order $n$ such that $\lambda_{\alpha}(G)\geq \lambda_{\alpha}(\mathcal{G}'_n)$ and $G\notin \mathcal{G}'_n$. Let $N$ be a sufficiently large integer such that Lemma \ref{L3} and Lemma \ref{L4} hold for $n> N$. Moreover, according to (\ref{4}), we assume $\lambda_{\alpha}(\mathcal{G}'_n)\geq (1-\varepsilon)\pi(\mathcal{F})n$ for $n> N$. Let $n_{0}=\left(\frac{Ne}{(1-\varepsilon)\pi(\mathcal{F})}\right)^{\frac{1}{\varepsilon^3}}$. Note that $n_0 > N$. Let $\mathbf{x}^{(n)}=\left(x^{(n)}_{1}, \dots, x^{(n)}_{n}\right)$ be a non-negative unit eigenvector of $A_{\alpha}(G)$ corresponding to $\lambda_{\alpha}(G)$.
Suppose $w^{(n)}\in V(G)$ such that $x_{w^{(n)}}=\min \left \{x^{(n)}_{1}, \dots, x^{(n)}_{n} \right\}$. Recall that $\lambda_{\alpha}(G)\geq \lambda_{\alpha}(\mathcal{G}'_n)$ and $\delta(G)\leq (\pi(\mathcal{F})-\varepsilon)n$. Then, by Lemma \ref{L3}, $x^2_{w^{(n)}}<\frac{1-\varepsilon^2}{n}$. For $n\geq n_0$, denote $G_{n}=G$ and $G_{n-1}=G_{n}-w^{(n)}$. Then, it follows from Lemma \ref{L4} that
\begin{equation*}
   \lambda_{\alpha}(G_{n-1})\geq \lambda_{\alpha}(G_n)\left(1-\frac{1-\varepsilon^3}{n-1}\right)
\end{equation*}
     and
\begin{equation*}
     \lambda_{\alpha}(G_{n-1})>\lambda_{\alpha}(\mathcal{G}'_{n-1}).
\end{equation*}

Since $\lambda_{\alpha}(G_{n-1})>\lambda_{\alpha}(\mathcal{G}'_{n-1})$, $G_{n-1}\notin \mathcal{G}'_{n-1}$. So $\delta(G_{n-1})\leq (\pi(\mathcal{F})-\varepsilon)(n-1)$. Let $\mathbf{x}^{(n-1)}=\left(x^{(n-1)}_{1}, \dots, x^{(n-1)}_{n-1}\right)$ be a non-negative unit eigenvector of $A_{\alpha}(G_{n-1})$ corresponding to $\lambda_{\alpha}(G_{n-1})$. Suppose $w^{(n-1)}\in V(G_{n-1})$ such that $x_{w^{(n-1)}}=\min \left \{x^{(n-1)}_{1}, \dots, x^{(n-1)}_{n-1} \right\}$. Denote $G_{n-2}=G_{n-1}-w^{(n-1)}$. Continue this operation, we can obtain a sequence of graphs $G_n$, $G_{n-1}$, $\dots$, which satisfies
\begin{equation*}
   \lambda_{\alpha}(G_{i-1})\geq \lambda_{\alpha}(G_i)\left(1-\frac{1-\varepsilon^3}{i-1}\right)
\end{equation*}
     and
\begin{equation*}
     \lambda_{\alpha}(G_{i-1})>\lambda_{\alpha}(\mathcal{G}'_{i-1})
\end{equation*}
for each $n\geq i> N$.

Then, for $i=N+1$, we have
\begin{eqnarray*}
    \lambda_{\alpha}(G_{N})
     &\geq& \lambda_{\alpha}(G_{N+1})\left(1-\frac{1-\varepsilon^3}{N}\right) \nonumber \\[2mm]
      &\geq& \lambda_{\alpha}(G_{n})\prod _{i=N+1}^n\left(1-\frac{1-\varepsilon^3}{i-1}\right) \nonumber \\[2mm]
       &\geq& \lambda_{\alpha}(G_{n})\mathrm{exp}\left(-\sum_{i=N+1}^{n}\left(\frac{1-\varepsilon^3}{i-1}+\left(\frac{1}{i-1}\right)^2\right)\right) \nonumber \\[2mm]
       &\geq & \lambda_{\alpha}(G_{n})\mathrm{exp}\left(-(1-\varepsilon^3)\ln \frac{n}{N-1}-1\right)  \nonumber \\[2mm]
       &\geq & (1-\varepsilon)\pi(\mathcal{F})n\left(\frac{n}{N-1}\right)^{-(1-\varepsilon^3)}e^{-1} \\[2mm]
        &\geq & (1-\varepsilon)\pi(\mathcal{F})n^{\varepsilon^3}e^{-1}\\[2mm]
       &\geq & (1-\varepsilon)\pi(\mathcal{F})n_{0}^{\varepsilon^3}e^{-1}=N,
\end{eqnarray*}
where the third inequality holds as Fact 1 and the fourth inequality holds as Fact 2. This contradicts the fact that $\lambda_{\alpha}(G_{N})\leq N-1$. The proof  is complete. \qed

\subsection{Proof of Theorem \ref{TT2}}

 Andr\'asfai, Erd\H{o}s,  and S\'os \cite{A} showed that $K_{r+1}$ is degree-stable with respect to the family of $r$-partite graphs. Erd\H{o}s and Simonovits \cite{E3}  extended the result to all color-critical graphs.

\begin{lem} [Erd\H{o}s and Simonovits \cite{E3}] \label{L5}
    Let $F$ be a color-critical graph with $\chi(F)=r+1\geq 3$. There is $n_0$ such that if $G$ is an $F$-free graph on $n\geq n_0$ vertices with $\delta(G)>\frac{3r-4}{3r-1}n$, then $G$ is $r$-partite.
\end{lem}

Next, we derive an upper bound for $\lambda_{\alpha}(T_{n,r})$. To this end, we first state the following lemma.

\begin{lem} \label{LL1}
     Let $\mathcal{G}$ be a hereditary and multiplicative graph family, and  let $r
     \geq 2$. If for any $0\leq \alpha\leq  1-\frac{1}{r}$, the inequality $\lambda_{\alpha}(\mathcal{G}_n)>\left(1-\frac{1}{r}\right)n-\left(1-\frac{1}{r}\right)$ holds, then
     $$ \lambda_{\alpha}(G)\leq \pi_{\alpha}(\mathcal{G})n$$
     for every $G\in \mathcal{G}_n$.
\end{lem}

\noindent\textbf{Proof.}
Let $H\in \mathcal{G}_n$   satisfy   $\lambda_{\alpha}(H)=\lambda_{\alpha}(\mathcal{G}_n)$. Since  $\mathcal{G}$  is multiplicative,  $H^p \in \mathcal{G}_{pn}$  for all   $p\geq 1$. Using the quotient matrix method, it can be verified that   $\lambda_{\alpha}(H^p)=p\lambda_{\alpha}(H)$. Consequently, $$\frac{\lambda_{\alpha}(\mathcal{G}_n)}{n}=\frac{\lambda_{\alpha}(H)}{n}=\frac{\lambda_{\alpha}(H^p)}{pn}\leq \frac{\lambda_{\alpha}(\mathcal{G}_{pn})}{pn}.$$ By Lemma \ref{T2} and the assumption  $\lambda_{\alpha}(\mathcal{G}_n)>\left(1-\frac{1}{r}\right)n-\left(1-\frac{1}{r}\right)$,  we have  $\lim_{p\to \infty} \frac{\lambda_{\alpha}(\mathcal{G}_{pn})}{pn}=\pi_{\alpha}(\mathcal{G})$. For any   $G\in \mathcal{G}_n$, we thus have $$\frac{\lambda_{\alpha}(G)}{n}\leq \frac{\lambda_{\alpha}(\mathcal{G}_n)}{n}\leq \pi_{\alpha}(\mathcal{G}), $$ which implies that $ \lambda_{\alpha}(G)\leq \pi_{\alpha}(\mathcal{G})n$. \qed

\begin{lem} \label{L6}
   For $r\geq 2$ and $0\leq \alpha\leq  1-\frac{1}{r}$, we have $\lambda_{\alpha}(T_{n,r})\leq \left(1-\frac{1}{r}\right)n$.
\end{lem}

\noindent\textbf{Proof.}
Let $\mathcal{G}$   denote the collection of all   $K_{r+1}$-free graphs. It is straightforward that $\mathcal{G}$   is both hereditary and multiplicative. Since   $T_{n,r}\in \mathcal{G}_n$, $\lambda_{\alpha}(\mathcal{G}_n)\geq \left(1-\frac{1}{r}\right)n-\frac{r}{4n}$. By Theorem \ref{TT1},   $\pi_{\alpha}(\mathcal{G})=1-\frac{1}{r}$, and thus Lemma \ref{LL1} yields $\lambda_{\alpha}(T_{n,r})\leq \left(1-\frac{1}{r}\right)n$. \qed

We are now ready to give the proof of Theorem \ref{TT2}.

\noindent\textbf{Proof of Theorem \ref{TT2}.}
Let $F$ be a color-critical graph with $\chi(F)=r+1\geq 3$, and let  $G$ be an arbitrary $F$-free $n$-vertex
graph. Recall that $\mathrm{ex}(n,F)=e(T_{n,r})=\left(1-\frac{1}{r}\right)\frac{n^2}{2}+O(1)$ and $\pi(F)=1-\frac{1}{r}$. Let $\varepsilon>0$ and $\sigma<\frac{\varepsilon^3}{7}$ be two sufficiently small constants. For sufficiently large $n$, we have
\begin{equation*}
    \begin{aligned}
       \left|\mathrm{ex}(n,F)-\mathrm{ex}(n-1,F)-\pi(F)n\right|=O(1)\leq \sigma n.
    \end{aligned}
\end{equation*}

Let $\mathcal{G}'_n$ be the set of all $F$-free graphs on $n$ vertices with minimum degree more  than $(\pi(F)-\varepsilon)n$. Note that $1-\frac{1}{r}-\varepsilon> \frac{3r-4}{3r-1}$ as $\varepsilon$ is sufficiently small. Then, for sufficiently large $n$, every graph in $\mathcal{G}'_n$ is an $r$-partite graph by Lemma \ref{L5}. Clearly, $T_{n,r}\in \mathcal{G}'_n$, which implies that $\lambda_{\alpha}(\mathcal{G}'_n)=\lambda_{\alpha}(T_{n,r})$ by Theorem \ref{TT0}. Combining with (\ref{4}) and Lemma \ref{L6}, we deduce that $\lambda_{\alpha}(\mathcal{G}'_n)=\lambda_{\alpha}(T_{n,r})=\left(1-\frac{1}{r}\right)n+o(1)$. Hence,
\begin{equation*}
    \begin{aligned}
       \left|\lambda_{\alpha}(\mathcal{G}'_n)-\frac{2\mathrm{ex}(n,F)}{n} \right|=o(1)\leq \sigma.
    \end{aligned}
\end{equation*}
Therefore, by Theorem \ref{TT3}, $\lambda_{\alpha}(G) \leq \lambda_{\alpha}(\mathcal{G}'_n)$ with equality if and only if $G=T_{n,r}$. This completes the proof. \qed


\begin{thebibliography}{99}

\bibitem{A}
 B. Andr\'asfai, P. Erd\H{o}s, V.T. S\'os, On the connection between chromatic number, maximal clique and minimum degree of a graph, Discrete Math. 8 (1974) 205-218.

 \bibitem{B}
 J. Byrne, D.N. Desai, M. Tait, A general theorem in spectral extremal graph theory, arXiv preprint arXiv: 2401.07266 (2024).

 \bibitem{C1}
 M. Chen, A. Liu, X. Zhang, On the $A_{\alpha}$-spectral radius of graphs
 without linear forests, Appl. Math. Comput. 450 (2023) 128005.

 \bibitem{C2}
 M. Chen, S. Li, Z. Li, Y. Yu, X. Zhang, An $A_{\alpha}$-
Spectral Erd\H{o}s-S\'os Theorem, Electron. J. Combin. 30(3) (2023) \#P3.34.

 \bibitem{E1}
 P. Erd\H{o}s, A.H. Stone, On the structure of linear graphs, Bull. Am. Math. Soc. 52 (1946) 1087-1091.

 \bibitem{E2}
 P. Erd\H{o}s, M. Simonovits, A limit theorem in graph theory, Studia Sci. Math. Hung. 1 (1966) 51-57.

 \bibitem{E3}
 P. Erd\H{o}s, M. Simonovits, On a valence problem in extremal graph theory, Discrete Math. 5 (1973) 323-334.

 \bibitem{F1}
 M.A.A. de Freitas, V. Nikiforov, L. Patuzzi, Maxima of the Q-index: forbidden 4-cycle and 5-cycle, Electron. J. Linear Algebra. 26 (2013) 905-916.

 \bibitem{H1}
 B. He, Y. Jin, X. Zhang, Sharp bounds for the signless Laplacian spectral radius in terms of clique number, Linear Algebra Appl. 438 (2013) 3851-3861.

 \bibitem{PK}
 P. Keevash, J. Lenz, D. Mubayi, Spectral extremal problems for hypergraphs, SIAM J. Discrete Math. 28(4) (2014) 1838-1854.

  \bibitem{Katona}
  G. Katona, T. Nemetz, M. Simonovits, On a problem of Tur\'an in the theory of graphs, Mat. Lapok 15 (1964) 228-328.

 \bibitem{L1}
 S. Li, Y. Yu, On $A_{\alpha}$ spectral extrema of graphs forbidding even cycles, Linear Algebra Appl. 668 (2023) 11-27.

 \bibitem{L2}
 S. Li, Y. Yu, H. Zhang, An $A_{\alpha}$-spectral Erd\H{o}s-P\'osa theorem, Discrete Math. 346(9) (2023) 113494.

 \bibitem{LY}
 Y. Li, W. Liu, L. Feng, A survey on spectral conditions for some extremal graph problems, Adv. Math. (China) 51(2) (2022) 193-258.

 \bibitem{N4}
 V. Nikiforov, Bounds on graph eigenvalues II, Linear Algebra Appl. 427 (2007) 183-189.

 \bibitem{N3}
 V. Nikiforov, Spectral saturation: inverting the spectral Tur\'an theorem, Electron. J. Combin. 16(1) (2009) R33.

 \bibitem{N2}
 V. Nikiforov, A spectral Erd\H{o}s-Stone-Bollob\'as theorem, Comb. Probab. Comput. 18 (2009) 455-458.

 \bibitem{N1}
  V. Nikiforov, Merging the $A$- and $Q$-spectral theories, Appl. Anal. Discrete Math. 11(1) (2017) 81-107.

 \bibitem{S1}
 M. Simonovits, A method for solving extremal problems in graph theory, stability problems, In Theory of Graphs (Proc. Colloq., Tihany, 1966) (1968) 279-319.

 \bibitem{T}
  P. Tur\'an, On an extremal problem in graph theory, Mat. Fiz. Lapok 48 (1941) 436-452.

 \bibitem{Y2}
 X. Yuan, Maxima of the Q-index: forbidden odd cycles, Linear Algebra Appl. 458 (2014) 207-216.

 \bibitem{Y1}
 X. Yuan, Z. Shao, On the maximal $\alpha$-spectral radius of graphs with given matching number, Linear Multilinear Algebra. 71(10) (2023) 1681-1690.

 \bibitem{Z2}
 J. Zheng, Y. Li, H. Li, The signless Laplacian spectral Tur\'an problems for color-critical graphs, arXiv preprint arXiv: 2504.07852 (2025).

 \bibitem{Z1}
 J. Zheng, H. Li, L. Su, A signless Laplacian spectral Erd\H{o}s-Stone-Simonovits theorem, Discrete Math. 349 (2026) 114665.
\end{thebibliography}
\end{document}